\documentclass[11pt]{article}
\usepackage{amsfonts,amsbsy,amsmath,amssymb}
\usepackage{graphics}

\usepackage[lmargin=2.6cm, rmargin=2.6cm,tmargin=3.0cm,bmargin=4.0cm]{geometry}

\numberwithin{equation}{section}

\newcommand{\iso}{\cong} 

\newcommand{\la}{\langle}
\newcommand{\ra}{\rangle}

\newcommand{\scaps}[1]{{\scshape #1}}
\newcommand{\bfm}{\textbf}
\newcommand{\mcal}{\mathcal}

\newcommand{\mbs}{\boldsymbol}

\newcommand{\G}{\mcal{G}}

\newcommand{\raf}[1]{(\ref{#1})}

\newcommand{\comp}{\hbox{$<\kern -3pt >$}}
\newcommand{\ncomp}
		{\;\hbox{\hbox{/}\kern -9.5pt \hbox{$<\kern -3pt >$}}}

\newcommand{\meet}
	       {\hbox{$\wedge \kern -5.75pt \raise 1.5pt \hbox{$.$}\,$}}
\newcommand{\Meet}
	     {\hbox{$\bigwedge \kern -8pt \raise 0.75pt \hbox{$.$}\:$}}
\newcommand{\ld}
	       {\hbox{$< \kern -6pt \raise 2pt \hbox{$.$}\,$}}

\newcommand{\nfrc}{\not \kern -5pt \frc}
\newcommand{\rest}{\vbox{\hbox{$\:\kern -2pt\mathbin{\vert\kern-3.1pt\lower-1pt
   \hbox{$\mathsurround=0pt\mathchar"0012$}\kern-4pt}\:$}}}

\newcommand{\empt}{\emptyset}

\newcommand{\pcntda}{\lower5pt\hbox{$\stackrel{\subset}{\neq}$}}
\newcommand{\pcntdb}{\lower5pt\hbox{$\stackrel{\supset}{\neq}$}}

\newcommand{\gD}{\Delta}

\newcommand{\II}{{\bf I\kern -1pt I}}

\newcommand{\surj}{\vbox{\hbox{$\longrightarrow $
                  \kern -22pt \hbox{\lower 2.5pt  \hbox{\tiny onto}}
                  \kern -16pt \hbox{\raise 5pt  \hbox{\tiny 1-1}}
                  \kern 3pt}}}
\newcommand{\uarrow}[2]{\vbox{\hbox{$\longrightarrow $
                  \kern -22pt \hbox{\lower 2.5pt  \hbox{\tiny $#2$}}
                  \kern -16pt \hbox{\raise 5pt  \hbox{\tiny $#1$}}
                  \kern 10pt }}}

\newcommand{\qed}{\kern 5pt\vrule height8.5pt width4.5pt depth0pt}

\newcommand{\sbd}[2]{#1\langle #2\rangle}
\newcommand{\modi}[1]{#1^{\mbox{\tiny{\rm mod}}}}
\newcommand{\mdl}[1]{\mbox{\rm mod $#1$}}
\newcommand{\nsp}{\mbox{A$\!^*$S}} 

\newcommand{\bd}{\begin{description}}
\newcommand{\ed}{\end{description}}

\newcommand{\ben}{\begin{enumerate}}
\newcommand{\een}{\end{enumerate}}

\newenvironment{theorem}
{
\refstepcounter{equation}
\ \\
\noindent
\begin{it}
\noindent
\bfm{Theorem~\theequation}
}
{\end{it}\ \\}

\newenvironment{lemma}
{
\refstepcounter{equation}
\ \\
\noindent
\begin{it}
\noindent
\bfm{Lemma~\theequation}
}
{\end{it}\ \\}

\newenvironment{prop}
{
\refstepcounter{equation}
\ \\
\noindent
\begin{it}
\noindent
\bfm{Proposition~\theequation}
}
{\end{it}\ \\}

\newenvironment{cor}
{
\refstepcounter{equation}
\ \\
\noindent
\begin{it}
\noindent
\bfm{Corollary~\theequation}
}
{\end{it}\ \\}

\newenvironment{obs}
{
\refstepcounter{equation}
\ \\
\noindent
\begin{it}
\noindent
\bfm{Observation~\theequation}
}
{\end{it}\ \\}

\newcommand{\sm}{\setminus}

\def\proof{\par\noindent{{\sl Proof.\ }}}
\def\QED{$\blacksquare$}
\def\inQED{$\square$}

\newcommand{\sref}[1]{\ref{#1}}

\makeatletter
\def\section{\@ifstar\unnumberedsection\numberedsection}
\def\numberedsection{\@ifnextchar[
  \numberedsectionwithtwoarguments\numberedsectionwithoneargument}
\def\unnumberedsection{\@ifnextchar[
  \unnumberedsectionwithtwoarguments\unnumberedsectionwithoneargument}
\def\numberedsectionwithoneargument#1{\numberedsectionwithtwoarguments[#1]{#1}}
\def\unnumberedsectionwithoneargument#1{\unnumberedsectionwithtwoarguments[#1]{#1}}
\def\numberedsectionwithtwoarguments[#1]#2{%
  \ifhmode\par\fi
  \removelastskip
  \vskip 3ex\goodbreak
  \refstepcounter{section}%
  \noindent
  \leavevmode
  \begingroup
  \bfseries
  \thesection.\ 
  #2.
  \endgroup
  \addcontentsline{toc}{section}{%
    \protect\numberline{\thesection}%
    #1}%
  }
\def\unnumberedsectionwithtwoarguments[#1]#2{%
  \ifhmode\par\fi
  \removelastskip
  \vskip 3ex\goodbreak
  \noindent
  \leavevmode
  \begingroup
  \bfseries
  #2.\quad
  \endgroup
  \addcontentsline{toc}{section}{%
    #1}%
  }
\makeatother

\makeatletter
\def\subsection{\@ifstar\unnumberedsubsection\numberedsubsection}
\def\numberedsubsection{\@ifnextchar[
  \numberedsubsectionwithtwoarguments\numberedsubsectionwithoneargument}
\def\unnumberedsubsection{\@ifnextchar[
  \unnumberedsubsectionwithtwoarguments\unnumberedsubsectionwithoneargument}
\def\numberedsubsectionwithoneargument#1{\numberedsubsectionwithtwoarguments[#1]{#1}}
\def\unnumberedsubsectionwithoneargument#1{\unnumberedsubsectionwithtwoarguments[#1]{#1}}
\def\numberedsubsectionwithtwoarguments[#1]#2{%
  \ifhmode\par\fi
  \removelastskip
  \vskip 3ex\goodbreak
  \refstepcounter{subsection}%
  \noindent
  \leavevmode
  \begingroup
  \bfseries
  \thesubsection.\ 
  #2.
  \endgroup
  \addcontentsline{toc}{subsection}{%
    \protect\numberline{\thesubsection}%
    #1}%
  }
\def\unnumberedsubsectionwithtwoarguments[#1]#2{%
  \ifhmode\par\fi
  \removelastskip
  \vskip 3ex\goodbreak
  \noindent
  \leavevmode
  \begingroup
  \bfseries
  #2.\quad
  \endgroup
  \addcontentsline{toc}{subsection}{%
    #1}%
  }
\makeatother

\begin{document}

\begin{center}
{\Large \scaps{Almost series-parallel graphs: structure and colorability}}\\ \ \\

{\small
Elad Aigner-Horev\\
horevel@cs.bgu.ac.il\\ \ \\

Department of Computer Science\\
Ben-Gurion University of the Negev\\ 
Beer Sheva, 84105, Israel}
\end{center}

{\small
\noindent
\bfm{Abstract.}
The series-parallel (SP) graphs are those containing no topological $K_{_4}$ and are considered trivial. We relax the prohibition distinguishing the SP graphs by forbidding only embeddings of $K_{_4}$ whose edges with both ends $3$-valent (skeleton hereafter) induce a graph isomorphic to certain prescribed subgraphs of $K_{_4}$. In particular, we describe the structure of the graphs containing no embedding of $K_{_4}$ whose skeleton is isomorphic to $P_{_3}$ or $P_{_4}$.
Such ``almost series-parallel graphs'' (ASP) still admit a concise description. 
Amongst other things, their description reveals that:
\begin{enumerate}

\item Essentially, the $3$-connected ASP graphs are those obtained from the $3$-connected cubic graphs by replacing each vertex with a triangle (e.g., the $3$-connected claw-free graphs). 

\item Except for $K_{_6}$, the ASP graphs are $5$-colorable in polynomial time.
Distinguishing between the $5$-chromatic and the $4$-colorable ASP graphs is $NP$-hard.

\item The ASP class is significantly richer than the SP class: $4$-vertex-colorability, $3$-edge-colorability, and Hamiltonicity are $NP$-hard for ASP graphs.
\end{enumerate}}

Our interest in such ASP graphs arises from  a previous paper of ours:``{\sl On the colorability of graphs with forbidden minors along paths and circuits}, Discrete Math. (to appear)''. \\

\noindent 
{\small \scaps{Keywords.} Topological $K_{_4}$, series-parallel graphs, graph coloring.}\\

\noindent
{\small \scaps{Preamble.} Whenever possible notation and terminology are that of~\cite{diestel}. Throughout, a graph is always simple, undirected, and finite. $G=(V,E)$ always denotes a graph. 

The skeleton of $G$, denoted $G^*$, is its subgraph induced by its edges with both ends at least $3$-valent. By $V^* \subseteq V$ we denote the set of vertices at least $3$-valent, and by $E^*$ the edges spanned by $V^*$; so that $G^*=G[E^*]$.  

A subdivision of a graph $H$ with its skeleton isomorphic to $S$ is denoted by $\sbd{H}{S}$. If $S= \emptyset$, then each edge is subdivided. 

By $W_{_r}$, $r \geq 2$, we denote the $r$-wheel consisting of an $r$-circuit, called the \emph{rim}, and an $r$-valent vertex, called the \emph{hub}, adjacent to every vertex of the rim through edges called \emph{spokes}. By $S_{_r}$ we denote a subdivided $W_{_r}$ with the spokes preserved. The hub and rim of an $S_{_r}$-subgraph are defied analogously. 

We call $G$ {\em virtually $3$-connected} if it is $2$-connected and for any two distinct $x,y\in V^*$ it contains three internally disjoint $xy$-paths.}

\section{Introduction} 
In this paper, we consider the structure of graphs distinguished by certain forbidden embeddings of $K_{_4}$. The forbidden embeddings are specified through constraints imposed on the skeleton of such an embedding. In particular, we describe the \emph{almost series-parallel} (ASP hereafter) graphs which are those containing no $\sbd{K_{_4}}{P_{_3}}$ or $\sbd{K_{_4}}{P_{_4}}$. 

Such supersets of the series-parallel (SP hereafter) class arise from~\cite{HOREV1} as follows. 
The starting point of~\cite{HOREV1} is the observation that graphs whose circuits induce subgraphs that are $K_{_r}$-minor free satisfy $\chi\leq f(r)$ for some $f$. 
It is then shown that graphs whose paths induce $K_{_{3,3}}$-minor free graphs are $6$-colorable; but the tightness of this result is unknown. Consequently, the following families were introduced.
Let $\G$ denote the class of graphs whose paths induce $K^-_{_{3,3}}$-minor free graphs and let $\G'$ denote the class of graphs whose paths induce $\{K_{_5},K^-_{_{3,3}}\}$-minor free graphs, where $K^-_{_{3,3}}$ is $K_{_{3,3}}$ with a single edge removed. Members of $\G'$ satisfy $\chi \leq 4$~\cite{HOREV1}. A query of~\cite{HOREV1} suggests that, in fact, except for $K_{_5}$ the members of $\G$ satisfy $\chi \leq 4$. 

It turns out that $\G\subseteq$ ASP (in fact $ \G \subseteq$ \nsp P; see below). To see this, observe that a topological $K^-_{_{3,3}}$ is in fact a $\sbd{K_{_4}}{S}$ with $\emptyset \subseteq S \subseteq C_{_4}$; and that a topological $K_{_4}$ is spanned by a path if and only if its skeleton is nonempty. Hence, every topological $K_{_4}$ in a member of $\G$ is a $\sbd{K_{_4}}{S}$ with $S$ either empty or containing a member of $\{K_{_3},K_{_{1,3}}\}$. 
Amongst other things, our results here regarding the colorability of ASP graphs settle the query of~\cite{HOREV1} to the affirmative. 

Through the notion of a {\sl skeleton} an alternative definition for classes such as $\G$ and $\G'$ that is more concise and less obscure is obtained. This is not the only merit of using such a ``language''. Such a ``language'' allows us access to certain tools and structures of graph theory that now become easier to apply. In our case, the use of graph connectivity as detailed in Section~\ref{decomposition} becomes useful and through which we reach the notion of {\sl virtual connectivity} (see Preamble).  
         
Such encounters with supersets of the SP graphs bring forth additional questions.
The SP graphs admit a simple structure~\cite{dir,duff} and are considered trivial as most classical graph theoretic problems are efficiently determined for these graphs~\cite{arn1,arn2,arn3}.
Limitations imposed on the skeleton of embeddings of $K_{_4}$ define 
supersets of the SP class. We may now think of an experiment which starts with the SP graphs and proceeds by relaxing the prohibition distinguishing the SP graphs as additional forms of embeddings of $K_{_4}$ become allowed; forbidden remain those embeddings whose skeleton satisfies certain constraints. We perform such an experiment with the somewhat obscure goal of ``finding where does graph theory begin''.  Eventually, we reach the ASP graphs. Of our experiment, we find the following hierarchy to be most interesting:
\begin{equation}\label{classes}
SP \subset \nsp P \subset ASP, 
\end{equation}
where a graph is \nsp P if it contains no $\sbd{K_{_4}}{P_{_3}}$, $\sbd{K_{_4}}{P_{_4}}$, or $\sbd{K_{_4}}{C_{_4}}$. 
Throughout \raf{classes}, we seek to observe:\\
(i) structural development; and\\
(ii) emergence of nontriviality in the form of NP-hard graph theoretic problems.\\ 

\noindent
\emph{Structural development.}  The SP class contains no $3$-connected members. The most notable structural development in (\ref{classes}) is the existence of a nontrivial class of $3$-connected ASP graphs; such admits a concise description (Theorem~\ref{A}). Actually, such a nontrivial class of $3$-connected graphs appears already at the \nsp P class.\\

\noindent
\emph{NP-hard problems.} We choose to track the computational complexity of the classical decision problems of: vertex colorability, edge colorability, and Hamiltonicity. All such decision problems can be solved efficiently for SP graphs~\cite{arn1,arn2,arn3}. For \nsp P graphs (and thus for ASP graphs), determining $3$-edge-colorability and Hamiltonicity is NP-hard. In particular, we observe the following phenomena:\\    
(i) for \nsp P graphs, $3$-vertex-coloring is hard while $4$-vertex-coloring can be solved efficiently; and\\ 
(ii) for ASP graphs, $4$-vertex-coloring is hard while $5$-vertex-coloring can be solved efficiently.\\
By this phenomena, we expect the class $ASP \sm \nsp P$ to be nontrivial.\\  
 
\noindent
\emph{Beyond ASP graphs.} If to extend the ASP graphs, then we recall that every graph of minimum degree $\geq 3$ contains a topological $K_{_4}$ whose skeleton contains a $P_{_4}$ (not necessarily isomorphic)~\cite{TT81,voss}. The description of the ASP graphs then suggests that perhaps most graphs satisfying $\delta \geq 3$ contain a $\sbd{K_{_4}}{P_{_4}}$.\\

\noindent
\scaps{Structure and colorability of $3$-connected ASP graphs.} 
A cubic graph is called \emph{truncated cubic} if each of its vertices belongs to precisely one triangle. Such are the graphs obtained from the cubic graphs by replacing each vertex with a triangle. Let $D_{_r}$ denote $K_{_{3,r}}+F$, where $F$ is a triangle spanned by the $3$-part of $K_{_{3,r}}$. 

The following theorem lists the $3$-connected ASP graphs together with their chromatic number. It asserts that the class of $3$-connected ASP graphs essentially coincides with the class of $3$-connected truncated cubic graphs; the latter are simply the $3$-connected cubic claw-free\footnote{A graph is claw-free if it does not contain $K_{_{1,3}}$ as an induced subgraph.} graphs.
 
\begin{theorem}\label{A}
A $3$-connected graph $G$ is ASP if and only if one of the following conditions is satisfied. 

\mbox{(\sref{A}.1)} $|V(G)|\le 6$ \mbox{\rm (obviously $2\le\chi(G)\le 6$)}, or
 
\mbox{(\sref{A}.2)} $K_{_{3,r}}\subseteq G\subseteq D_{_r}$, $r\ge 4$ \mbox{\rm (so 
that $2\le\chi(G)\le 4$)}, or

\mbox{(\sref{A}.3)} $G\cong W_{_r}$, $r\ge 6$ \mbox{\rm (so that $3\le\chi(G)\le 4$)}, or
 
\mbox{(\sref{A}.4)} $G$ is truncated cubic \mbox{\rm ($\chi(G)=3$, by Brooks' 
theorem)}.
\end{theorem}

By (\sref{A}.4), the problems of Hamiltonicity and $3$-edge-colorability are $NP$-hard for the class of $3$-connected ASP graphs~\cite{gj}. \\

\noindent
\scaps{Structure of $2$-connected ASP graphs.} The $2$-connected ASP graphs may be described in two steps. First, we decompose a $2$-connected ASP graph into virtually $3$-connected ASP graphs (using the standard decomposition of $2$-connected graphs into their $3$-connected components; see Section~\ref{decomposition}). Second, we characterize the virtually $3$-connected ASP graphs in terms of their skeleton components as specified in Theorem~\ref{B}.

We require the following notation and terminology. 
Except for Section~\ref{color-ASP}, throughout a graph $G=(V,E)$ is $2$-connected. A path in $G$ with its ends in $V^*$ and its interior vetices in $V\sm V^*$ is called a {\em thread}. The end-pair of a thread will be called a {\em window}.
Here, a graph is assumed to be {\sl simple}, in the sense that it has no parallel edges and no parallel threads (``parallel'' means having the same ends; a thread parallel to a skeleton edge is allowed). When the difference is irrelevant, the term {\em link} will refer to both edge and thread. The ends of a link are considered as {\sl neighbors}; by $N(v)$ we denote the set of so understood neighbors of a vertex $v$.

For $v \in V$, let $d^*(v) = d_{G^*}(v)$. For $v \in V^*$, introduce two subsets of neighbors, not necessarily disjoint:\\
$N^1(v)\colon=\{u\in V\colon\;uv\in E^*\}$, so that $d^*(v)= d_{G^*}(v)=|N^1(v)|$, and\\ 
$N^2(v)\colon=\{u\in V^*\colon\;\{u,v\}\;\;\mbox{\rm is a window}\}$.\\
Clearly, $N(v)=N^1(v)\cup N^2(v)$.

A subdivision of a simple $3$-connected graph is clearly virtually $3$-connected, but not only:

\begin{obs}\label{virt}
A graph with at least four vertices is virtually $3$-connected if and only if it is obtainable 
from a simple $3$-connected graph $(V,E)$ by choosing sets $F\subseteq W\subseteq E$, 
attaching a thread in parallel to each member of $W$, and removing the set $F$ after
that. The skeleton of the resulting graph is $E\sm F$, and the set of windows is the set of pairs defined by $W$.
\end{obs} 

The following two sequences of virtually $3$-connected ASP graphs are obtained by attaching threads in parallel to certain edges of graphs from Theorem~\sref{A}:

\noindent
$\modi{W_{_r}}\colon=W_{_r}+r$ threads attached to the rim edges;\\
$\modi{D_{_r}}\colon=D_{_r}+3$ threads attached to edges of the triangle spanned by the part of size $3$ in $K_{_{3,r}}$.

Finally, by {\sl fishpond} we mean a virtually $3$-connected graph whose skeleton components satisfy the following definition.\\

\noindent
{\sc Definition A.}\\
A virtually $3$-connected graph $G$ is a {\em fishpond} if $\gD(G^*)\le 4$, and 
\begin{enumerate}
\item [(A.1)] each vertex with $d^*(v)\ge 2$ belongs to a skeleton triangle;
\item [(A.2)] in a skeleton triangle, no edge has a parallel thread attached to its ends, and at least two of its vertices $v$ have $|N(v)|=3$;
\item [(A.3)] each component $J$ of $G^*$ satisfies one of the following conditions:

$\hspace{0.3cm}$ (A.3.1) $|V(J)|\le 3$;

$\hspace{0.3cm}$ (A.3.2) $\gD(J)=3$;

$\hspace{0.3cm}$ (A.3.3) $J$ is a union of two edge-disjoint triangles with one vertex  in common.
\end{enumerate}

A total subdivision of a simple $3$-connected graph is
a fishpond (with no fish). Since a fishpond is virtually $3$-connected, if it has at least $5$ skeleton vertices, then it contains no $K_{_4}^-$. For suppose $J \iso K^-_{_4}$ is contained in a fishpond, and suppose $t$ is a skeleton vertex external to $J$. A $(t,J)$-fan\footnote{See definition of {\sl fans} in Section~\ref{pre}.} of order $3$ clearly has one of its ends at a $3$-valent vertex of $J$, say $u$, implying that $|N(u)| > 3$; contradicting (A.3.2).  

Since a $3$-connected graph is a virtually-$3$-connected graph coinciding with its skeleton, a $3$-connected fishpond is a truncated cubic graph. Indeed, a $3$-connected fishpond has a single skeleton component which is itself so that $d^*(v) \geq 3$ for each $v \in V(G)$, by $3$-connectivity. (A.3.2) then implies that the graph is cubic. By (A.1) each vertex belongs to a triangle. But since $G$ is cubic and contains no $K^-_{_4}$ (if $|V^*|\geq 5$), each 
vertex belongs to precisely one triangle. 

A skeleton component $J$ of a fishpond has one of the following forms. 
In case (A.3.1), $|V(J)|\le 3$ clearly implies that $J\cong K_{_i}$, $i=1,2,3$.
In case (A.3.2), $J$ is a connected subgraph of a truncated cubic graph, with 
$|V(J)|\ge 4$ (as there is no $K^-_{_4}$). In case (A.3.3), each $2$-valent vertex of $J$ has $|N(v)|=3$. For the $4$-valent vertex $v$ in $J$, we have $|N(v)|=4$. This is since, by (A.2), no $(t,J)$-fan of order $3$, where $t \in V^* \sm V(J)$, can end in the $4$-valent vertex and a $2$-valent vertex of $J$; if $|N(v)| \geq 5$, such a fan exists. 

\begin{theorem}\label{B}
A virtually $3$-connected graph $G=(V,E)$ with $|V^*| > 6$ is ASP if and only if it satisfies one of the following conditions:

\mbox{(\sref{B}.1)}
$K_{_{3,r}}\subseteq G\subseteq\modi{D_{_r}}$, $r\ge 4$; 

\mbox{(\sref{B}.2)}
$S_{_r}\subseteq G\subseteq\modi{W_{_r}}$, $r\ge 6$; 

\mbox{(\sref{B}.3)}
$G$ is a fishpond.
\end{theorem}

Together with \raf{rcpt} (the decomposition), Theorem~\ref{B} characterizes the ASP graphs, and 
imply that identification of ASP graphs can be done efficiently. 
Since a $3$-connected fishpond is a truncated cubic graph, the list of $3$-connected ASP graphs (Theorem~\sref{A}) is immediately extractable from Theorem~\sref{B}, because a $3$-connected graph is just a virtually $3$-connected graph coinciding with its skeleton. \\

\noindent
\scaps{Colorability of ASP graphs.}

\begin{theorem}\label{C}
Apart from $K_{_6}$, the ASP graphs are $5$-vertex-colorable in polynomial time.
\end{theorem}

There are $5$-chromatic ASP graphs not containing $K_{_5}$ as a subgraph (see Section~\ref{color-ASP}). Characterizing the $5$-chromatic ASP graphs would be a difficult problem as indicated by the following.   

\begin{theorem}\label{D}
Deciding $4$-colorability of ASP graphs is NP-complete. 
\end{theorem}

\noindent
\scaps{Colorability of \nsp P graphs.}   

\begin{prop}\label{E}
Apart from $K_{_5}$, the \nsp P graphs are $4$-vertex-colorable in polynomial time.
\end{prop}

\noindent
There are planar $4$-chromatic \nsp P graphs not containing $K_{_4}$ as a subgraph (see Section~\ref{color-ASP}).

\begin{prop}\label{F}
Deciding $3$-colorability of planar \nsp P graphs is NP-complete. 
\end{prop}

\noindent
\scaps{An application.} We have seen that $\G \subseteq \nsp P$. Proposition~\ref{E} then settles the query of~\cite{HOREV1} mentioned above (note that $K_{_4} \in \G$).

\begin{cor}
\label{k33-path-2}
Apart from $K_{_5}$, the members of $\G$ are $4$-colorable.
\end{cor}

\noindent
\scaps{Outline of the paper.} 
A virtual $3$-connected ASP graph with $|V^*| \leq 4$ is trivial. Hence, we always assume that 
\begin{equation}\label{>4}
|V^*| > 4
\end{equation}
To characterize the virtually $3$-connected ASP graphs we proceed in two steps. First, in Section~\ref{structure-asp'}, we describe the virtually $3$-connected ASP graphs containing a $\sbd{K_{_4}}{C_{_4}}$-subgraph. 
Second, in Section~\ref{asp-structure}, we characterize the virtually $3$-connected \nsp P graphs. These two characterizations, are then combined in Section~\ref{proof-B} to yield a proof of Theorem~\ref{B}. Finally, Colorability of ASP graphs is discussed in Section~\ref{color-ASP}.

\section{Preliminaries}\label{pre}
Let $H \subseteq G$ be a subgraph of $G$. An $H$-\emph{ear} is a path internally disjoint of $H$;
if its ends are $x,y \in V(H)$, then such an ear is also called an $(x,y)$-{\sl ear} of $H$.

A link (and also its end-pair) in a virtually $3$-connected ASP graph will be called 
{\em red} if attaching a parallel link (of the opposite kind) creates a forbidden 
subgraph. An ASP graph with all the members of $E^*\sm W$ red will be called {\sl maximal}, where $W$ is as in Observation~\ref{virt}.

Let $v\in V$ and $H$ be a subset of $V\sm\{v\}$ or a subgraph of $G-v$. A union of $(v,H)$-paths disjoint in $G-v$ and having no common vertex with $H$, except the ends, will be called a $(v,H)$-{\em fan}. The {\em ends} of a fan are its 
common vertices with $H$. 

By the Menger theorem, a $(v,H)$-fan with $k$ ends (or 
{\sl of order} $k$) exists if and only if $G$ has no $(v,H)$-disconnector of cardinality $<k$. Moreover, if in such a case there is a $(v,H)$-fan $F'$ with the end-set $X$, $|X|<k$, then there also exists a $(v,H)$-fan $F$ of order $k$ whose end-set contains $X$. Such an $F$ may be chosen so as to contain those $(v,H)$-paths of $F'$ which are links.  
Throughout this chapter, the order of a fan is always $3$. When $G$ is virtually $3$-connected, a $(v,H)$-fan of order $3$ exists for each $v\in V^*\sm H$, provided $|V(H)\cap V^*|\ge 3$.

For vertices $x$ and $y$ of a graph $H$, we denote by $[xHy]$ an $xy$-path in $H$ when
such a path is unique (as in a tree), or may be chosen arbitrarily, or is clear from 
the context (e.g., when three or more vertices are specified in a circuit $C$, $[xCy]$ is the minimal segment containing exactly two of those vertices, namely $x$ and $y$). Such a path may also be considered as semi-open or open, denoted by $[xHy)$, $(xHy]$, and
$(xHy)$; this means excluding one end or both: for instance, $[xHy)$ is an alternating sequence of edges and vertices with $x$ on one end and an edge incident with $y$ on the other.
Also, if $P$ is a $uv$-path, then we write $int P$ to denote $V(P) \sm \{u,v\}$. 

We conclude this section with the following easy observation. 

\begin{obs}\label{whi}
A $2$-connected graph with no circuit traversing vertices $x,y,z$ contains a subdivision 
of $K_{_{3,2}}$ with $\{x,y,z\}$ as the larger part.
\end{obs}

\section{Receptacles}\label{decomposition}
Through the standard decomposition of a $2$-connected graph into its $3$-connected components~\cite[Section 9.4]{bm}, a $2$-connected ASP graph can be decomposed into virtually $3$-connected ASP graphs as follows.

Given a graph $G$, a maximal virtually $3$-connected subgraph of $G$ is called
a {\em receptacle}. 
The skeleton of a receptacle is the set of edges spanned by a maximal set $X\subseteq V(G)$ such that $X\sm\{u,v\}$ is not disconnected in $G-\{u,v\}$, for whatever $u,v\in V(G)$, and its windows are the pairs $\{u,v\}\subseteq X$ disconnecting $G$. Thus, up to the lengths and exact location of the threads outside $X$, the receptacles are uniquely defined by $G$. 
The intersection of two receptacles is an induced subgraph of $G$ of order $\le 2$; a window may belong to more than two receptacles. 

Thus, the receptacle of $G$ are its $3$-connected components with their ``marked edges'' (see ~\cite[Section 9.4]{bm}) subdivided. Consequently, the receptacles of a given $2$-connected graph are available through the procedure described in~\cite[Section 9.4]{bm}.

Trivially, an embedding of a $3$-connected graph ($K_{_4}$ in our case) into a $2$-connected graph $G$ has all its skeleton vertices contained in exactly one receptacle. Therefore,

\begin{equation}\label{rcpt}
\mbox{\rm a graph is ASP if and only if each of its receptacles is ASP.}
\end{equation}

\section{ASP graphs with a $\sbd{K_{_4}}{C_{_4}}$-subgraph}\label{structure-asp'}
In this section, we prove the following

\begin{prop}\label{kc4}
A virtually-$3$-connected ASP graph $G=(V,E)$ with $|V^*| \geq 7$ contains $\sbd{K_{_4}}{C_{_4}}$ if and only if $K_{_{3,r}}\subseteq G\subseteq\modi{D_{_r}}$, $r\ge 4$.
\end{prop}

\noindent
\scaps{Preliminaries.} Let $G$ be an ASP graph and $C$ be a $4$-circuit in $G$, with the vertices $u_{_i}\in V^*$, $i=0,1,2,3$ $\mdl{4}$, so ordered along $C$. 
Let $J$ be an embedding of $K_{_4}$ in $G$, with the skeleton $C$; so that $J \iso \sbd{K_{_4}}{C_{_4}}$. For $i=1,2$, the open $(u_{_{i-1}},u_{_{i+1}})$-path in $J-E(C)$, clearly contains $\geq 2$ edges (since the skeleton of $J$ is precisely $C$), will be 
called a \emph{diagonal} and denoted by $L_{_i}$. Recall that we assume $|V^*|>4$; consequently,

\begin{equation}\label{thds}
\mbox{\rm at least one diagonal of $J$ meets $V^*$,}
\end{equation}
To see \raf{thds}, note that otherwise $V^*\subseteq J$. Indeed, suppose, to 
the contrary, that there is an $s\in V^*\sm J$, and construct $J'\colon=J\cup F$ where $F$ is an $(s,J)$-fan. Let the ends of $F$ be $u_{_i}$, $i=0,1,2$. Then, $(sFu_{_0})$ and $(sFu_{_2})$ are edges, the latter because otherwise $J'-u_{_0}u_{_1}-L_{_1}$ is 
$\sbd{K_{_4}}{S}$ with $S\cong P_{_3}$ or $P_{_4}$, and the former by the symmetry. But now $J'-\{su_{_0},u_{_1}u_{_2}\}$ is $\sbd{K_{_4}}{P_{_4}}$, contradiction. 

In fact, we observe something more systematic. 

\begin{obs}\label{svj}
Let $s\in V^*\sm C$ and $F$ be an $(s,C)$-fan with the ends $u_{_{i-1}}$, $u_{_i}$ and
$u_{_{i+1}}$. If $(sFu_{_i})$ is not an edge, then the other two members of $F$ 
are edges. 
\end{obs}

\noindent
\scaps{Preparations towards the ``IF'' direction of Proposition~\ref{kc4}.} 
The ``if'' direction of our proof of Proposition~\ref{kc4} proceeds by distinguishing two cases. If $G$ is as in the premise of Proposition~\ref{kc4}, then it has a subgraph $J$ as above. Either there is a $J$-ear linking the interior of its diagonals or there is not. In the former case, we have that $G$ contains a topological $K_{_{3,3}}$ with its skeleton containing $C$; we shall see that this implies that $K_{_{3,r}}\subseteq G\subseteq\modi{D_{_r}}$, for $r\ge 4$ (see Lemma~\ref{Dr}). In the latter case, it turns out that $|V^*| \leq 6$ (see Lemma~\ref{k23}).  

\begin{lemma}\label{k23}
Let $G$ be an ASP graph, and $C$ be a $4$-circuit in $G$ with $V(C)\subseteq V^*$. 
Suppose that $G$ contains a subdivision of $K_{_4}$ with skeleton $C$ and no subdivision of $K_{_{3,3}}$ whose skeleton contains $C$. Then, $|V^*| \leq 6$. 
\end{lemma}

\noindent
{\sl Proof.}
Let $J$, $C$, $L_{_1}$, and $L_{_2}$ be as above, and let $s\in L_{_2}\cap V^*$, by 
\raf{thds}. 

\noindent
(I) We show that $C\cup L_{_2}\cong K_{_{2,3}}$. Construct $J'\colon=J-L_{_2}\cup F$, where $F$ 
is an $(s,J-L_{_2})$-fan with two of its ends coinciding with the ends of $L_{_2}$ (i.e., $u_{_1}$ and $u_{_3}$). 
Since $C$ is contained by no subdivision of $K_{_{3,3}}$, let the third end of $F$ be
$u_{_2}$, by symmetry. Then, $(sFu_{_i})$, $i=1$ and $3$, are edges, for otherwise $J'-u_{_1}u_{_2}$ or $J'-u_{_2}u_{_3}$ is non-ASP. \\

\noindent
(II) We show that $P\colon=(sFu_{_2})$ is a link (edge or thread). Suppose, to the contrary,
that there is a $t\in P\cap V^*$, and construct $J''\colon=J'-P\cup F'$ where $F'$ is
a $(t,J''-P)$-fan with two of its ends coinciding with the ends of $P$ (i.e., $s$ and $u_{_2}$). Denote the third 
end of $F'$ by $v$. We have $v\not\in L_{_1}\cup\{u_{_0}\}$, for otherwise $J''$ contains
$\sbd{K_{_{3,3}}}{S}$ with $C\subseteq S$. In addition, $v\ne u_{_1},u_{_3}$; for if,
say, $v=u_{_1}$, then $J''-\{su_{_1},u_{_1}u_{_2}\}\cong\sbd{K_{_4}}{P_{_i}}$,
$i=3$ or $4$. Thus, $P\cap V^*=\empt$, as required.\\

\noindent
(III) Suppose that there is a $t\in int L_{_1}\cap V^*$; we show that $|V^*|\leq 6$. Note that if $P$ is a thread, then $int L_{_1} \cap V^* =\emptyset$ by symmetry. Thus, by the assumption that $t\in int L_{_1}\cap V^*$, we may assume that $P$ is an edge. We may proceed as above, namely, construct 
$J''\colon=J'-L_{_1}\cup F'$ where $F'$ is a $(t,J'-L_{_1})$-fan with two of its ends coinciding with the ends of $L_{_1}$. Let $v$ be the third end of $F'$. Since $G$ is assumed to contain no $\sbd{K_{_{3,3}}}{S}$ with $C\subseteq S$, we have $v\ne s$. Thus, $J''$ has two versions: $A$, with $v=u_{_1}$, and $B$, with $v=u_{_3}$. We have $A-su_{_2}\cong B-su_{_2}$ and $A-\{u_{_1}u_{_2},u_{_2}u_{_3}\}\cong B-\{su_{_1},u_{_2}u_{_3}\}$. 

Since $A-\{su_{_2},u_{_0}u_{_1}\}$ and $A-\{su_{_2},u_{_1}u_{_2}\}$ are ASP, the paths $(tF'u_{_0})$ and $(tF'u_{_2})$ are edges, in both $A$ and $B$. Since 
$A-\{u_{_1}u_{_2},u_{_2}u_{_3}\}$ is ASP, $(tF'u_{_1})$ should be an edge, in both $A$ and $B$. Thus, $A$ and $B$ coincide each with its skeleton; implying that $|V^*| \leq 6$. \\

\noindent
(IV) To conclude, Suppose then that $int L_{_1} \cap V^* = \emptyset$. 
Let, to the contrary, $t \in V^* \sm V(J)$. By the assumption of this case and (II), $P$ is a link and so a $(t,J)$-fan $F$ has its ends contained in $\{u_{_0},u_{_1},u_{_2},u_{_3},s\}$. 
It is routine to check that $J \cup F$ contains a forbidden configuration; implying that $t$ does not exist and so $|V^*| \leq 6$ as required. \QED \\

Prior to Lemma~\ref{Dr}, we note that Observation~\ref{svj} implies that  
\begin{equation}\label{*}
\mbox{\rm an ASP graph $H\iso \sbd{K_{_{3,3}}}{S}$ with $C_{_4} \subseteq S$ satisfies $H \iso K_{_{3,3}}$.}
\end{equation}

We also observe that 
\begin{equation}\label{***}
\mbox{the only maximal ASP graph with the skeleton containing $K_{_{3,r}}$, $r\ge 4$, is $\modi{D_{_r}}$.}
\end{equation}
To see \raf{***}, note that a thread cannot be appended in parallel to an edge of a $K_{_{3,r}}$, $r \geq 3$, by \raf{*}. Next, no two vertices in the part of cardinality $r \geq 4$ can be connected by a link. Indeed, let $R \iso K_{3,4}$ and let $z$ and $w$ be vertices on the part of cardinality $4$. Let $y \in V(R)$ be a vertex on the part of cardinality $3$. It follows that a $(z,w)$-link (once added) and the edge $zy$ can be used to form a $\sbd{K_{3,3}}{S}$, $C_{_4} \subseteq S$, that is not isomorphic to $K_{3,3}$; contradicting \raf{*}. 

\begin{lemma}\label{Dr}
Let $G=(V,E)$, $|V^*| \geq 7$, be a virtually-$3$-connected ASP graph, and let $C$ be a $4$-circuit in $G$ with $V(C)\subseteq V^*$. 
If $G$ contains a $\sbd{K_{_4}}{S}$ with $S = C$, then $K_{_{3,r}}\subseteq G\subseteq\modi{D_{_r}}$, $r\ge 4$.
\end{lemma}

\proof
Let $J,C,L_{_1}$, and $L_{_2}$ be as above. We may assume that there is a topological $K_{_{3,3}}$ with its skeleton containing $C$, for otherwise $|V^*| \leq 6$ , by Lemma~\ref{k23}. Such a subdivision is isomorphic to $K_{_{3,3}}$, by \raf{*}.  
Let $K_{_{3,3}} \iso H \subset G$ and put $V(H) = \{x_{_i}\}_{i=1}^6$ such that the bipartition of $H$ is $A = \{x_{_1},x_{_2},x_{_3}\}$, $B = \{x_{_4},x_{_5},x_{_6}\}$. 
Let $X = V^* \sm V(H^*)$. \\

\noindent
(I) If $x \in X$ and $F$ is an $(x,H)$-fan, then the ends of $F$ are contained all in $A$ or all in $B$. To see this, assume, to the contrary, that the ends of $F$ are $x_1,x_2$, and $x_4$. The $4$-circuit $\la x_3,x_6,x_2,x_5 \ra$, $x$, and the subgraph $F+x_{_1}x_{_6}+x_{_4}x_{_3}$ contradict Observation~\ref{svj}.\\

\noindent
(II) Let $x$ and $F$ be as in (I). Then, the members of $F$ are edges. 
Indeed, by (I), the ends of $F$ are contained in, say, $A$. Apply Observation~\ref{svj} to $F$ and the circuits induced by 
$\{x_{_1},x_{_4},x_{_2},x_{_5}\}$ and $\{x_{_2},x_{_5},x_{_3},x_{_6}\}$.\\

\noindent
(III) Let $H' \iso K_{_{3,4}}$ with bipartition $A',B'$, $|A'|=3$. 
Suppose $V(H') \subseteq V^*$ and let $x \in V^* \sm V(H'^*)$. 
Then, an $(x,H')$-fan has all its ends contained in $A'$. 
To see this, assume to the contrary, that this is not the case, and thus, by (II), 
the ends of $F$ are contained in $B'$. Let $C'$ be a $4$-circuit in $H'$ 
containing a vertex in $B'$ that is not an end of $F$. 
The paths in $F$ can be extended through edges of $H'$ as to end in $C'$ and
to contradict Observation~\ref{svj}.\\

Claims (I),(II), and (III) imply that $G^*\iso K_{_{3,r}}$, where $r = 3+|X| \geq 4$. Claim (II), also asserts that no thread of $G$ is parallel to an edge of $G^*$. The claim now follows from \raf{***}.     
\QED \\

We are now ready to prove Proposition~\ref{kc4}.

\noindent
\bfm{Proof of Proposition~\ref{kc4}.} ``Only if'' is a routine check. To show ``if'', let $G$ be a virtually $3$-connected ASP graph containing $J,C$, and $L_{_i}$, $i=1,2$, as defined above. Either there is a $J$-ear linking the interiors of $L_{_1}$ and $L_{_2}$ or there is not. In the former case, $K_{_{3,r}}\subseteq G\subseteq\modi{D_{_r}}$, $r\ge 4$,
by Lemma~\ref{Dr}. In the latter case, $|V^*| \leq 6$, by Lemma~\ref{k23}, contradicting the assumption that $|V^*|\geq 7$. 
\QED 

\section{Structures in \nsp P graphs}\label{asp-structure}
In the previous section, structure of ASP graphs containing a $\sbd{K_{_4}}{C_{_4}}$ was considered (see Proposition~\sref{kc4}). It remains to describe the \nsp P graphs. In this section, we prepare towards the characterization of the \nsp P graphs which is (eventually) presented in our proof of Theorem~\ref{B} (see Section~\ref{proof-B}). Here, we study local structures surrounding a 
skeleton vertex (Section~\ref{vertex}) and a skeleton triangle (Section~\ref{triangle}).
These structures will be used in our proof of Theorem~\ref{B}.

\subsection{Preliminaries}

\begin{lemma}\label{xvy}
A virtually $3$-connected graph $G=(V,E)$ is non-\mbox{\rm \nsp P} if and only if it has a skeleton 
path $P=[x,v,y]$ and a circuit $C$ traversing $x$ and $y$ but neither $v$ nor an
$xy$-edge, such that

\begin{equation}\label{DN}
\left(N^2(v)\sm\{x,y\}\right)\cup\left(N^1(v)\sm C\right)\ne\empt.
\end{equation}
\end{lemma}

\noindent
{\sl Proof}
``Only if'' is straightforward: if $J$ is a forbidden $\sbd{K_{_4}}{S}$-subgraph with $S$ containing $P$, then \raf{DN} holds for $C\colon=J-v$. 

Suppose now that the ``if'' assertion is false, and let $G$ be an \nsp P graph containing $P=[x,v,y]$ and $C$ as specified such that \raf{DN} is satisfied. Then, clearly $N^2(v)$ does not meet $C-\{x,y\}$, and \raf{DN} implies that there is $z\in N(v)\sm C$. Now,
\begin{equation}\label{*'}
\mbox{\rm any $(z,C\cup P)$-fan $F$ in $G$ containing the $(v,z)$-link is incident with
$x$ and $y$,}
\end{equation}
for if such an $F$ terminates in $a\in C-\{x,y\}$, then $C\cup P\cup(aFv)$ is non-\nsp P.
Fix $F$ and put $J\colon=C\cup P\cup F$. The $(z,v)$-link should be an edge, for 
otherwise removing from $J$ one of the open $(x,y)$-segments of $C$ yields a non-\nsp P 
subdivision of $K_{_4}$.

Since $G$ is virtually $3$-connected, $J$ has an ear, say $Q$, connecting $C-\{x,y\}$ 
and $F-\{x,y\}$. By (\ref{*'}), the end of $Q$ in $F$ is just $v$; let $u$ denote its
other end, in $C$. Then, $(zFx)$ and $(zFy)$ are edges, the latter because otherwise 
$J\cup Q-vx-(uCy)$ is non-\nsp P and the former by symmetry. But now 
$J\cup Q-\{xv,yv\}$ is non-\nsp P, contradiction.
\QED
\\

The following particular cases are often encountered: 

\begin{equation}\label{spo}
\mbox{\rm for $r\ge 4$, the spokes of an $S_{_r}$-subgraph are red;}
\end{equation}

\begin{equation}\label{k33}
\mbox{\rm a subdivision of $K_{_{3,3}}$ with the skeleton containing $P_{_3}$ is 
non-\nsp P,}
\end{equation}
because $K_{_{3,3}}$ contains a union of a $4$-circuit $C$ and a $(v,C)$-fan with
the paths of lengths $1,1,2$.

\subsection{Structure surrounding a skeleton vertex}\label{vertex} 
The starting point of this section (and its main message) is that a vertex with $d^* \geq 2$ is the hub of an $S_{_{d^*}}$-subgraph. 

A vertex $v$ of a $2$-connected graph $G$ satisfying $d(v)>2$ and $\kappa(G-v) \geq 2$
is called \emph{regular}. Trivially, a vertex $v$ is regular if and only if $N^2(v)=\empt$.
Thus, a vertex of a virtually $3$-connected graph is regular if and only if it is incident with no thread.

\begin{lemma}\label{rd*}
Let $G=(V,E)$ be virtually $3$-connected and \mbox{\rm \nsp P}, and let $v\in V^*$ have 
$d^*(v)\ge 2$. Then,

\mbox{(\sref{rd*}.1)} $v$ is the hub of an $S_{_r}$-subgraph with $r=d^*(v)$; 

\mbox{(\sref{rd*}.2)} when $d^*(v)=2$, such a subgraph may be chosen so as to
contain a vertex from \indent \indent \indent \ $N^2(v)\sm N^1(v)$;

\mbox{(\sref{rd*}.3)} if $d^*(v)\ge 3$, then $N^2(v)\subseteq N^1(v)$;

\mbox{(\sref{rd*}.4)} if $d^*(v)\ge 4$, then $v$ is regular and $N^2(v) = \emptyset$.
\end{lemma}

\noindent
{\sl Proof.} 
Let us note that (\sref{rd*}.4) is implied by (\sref{rd*}.3) and \raf{spo};
so it will be sufficient to prove the first three terms of the lemma. 

To see (\sref{rd*}.2), choose $X\subseteq N(v)$ with $|X|=3$ and $|X\cap N^1(v)|$ 
as large as possible. Then, $G-v$ has a circuit containing $X$. 
Indeed, otherwise $G-v$ contains a subdivision of $K_{_{3,2}}$ with $X$ as the larger 
part, by Observation~\sref{whi}. Together with $\{v\}$ and the $(v,X)$-links, it forms
a subdivision of $K_{_{3,3}}$ whose skeleton contains $P_{_3}$; contradiction to \raf{k33}.

To prove (\sref{rd*}.1) and (\sref{rd*}.3), we may assume that $d^*(v)\ge 3$ (as (\sref{rd*}.1) is true for $d^*(v) =2$ by the proof of (\sref{rd*}.2) above). Let $C$ be 
a circuit of $G-v$ with $X\colon=C\cap N^1(v)$ maximal. We are to show that $X=N^1(v)$ 
and that $N^2(v)\sm N^1(v)=\empt$. Put $r\colon=|X|$, and let $x_{_i}$, $i=0,\dots,r-1$ 
$\mdl{r}$, be a circular ordering of $X$ along $C$. Let $J\cong S_{_r}$ be the union of
$C$, $\{v\}$ and the set of $(v,X)$-edges. As in the previous case, we may assume $|X|\ge 3$. Hence, no vertex from $N^2(v)\sm N^1(v)$ belongs to $C$, by (\ref{spo}). If there is a $u\in N(v)\sm C$, 
consider a $(u,J)$-fan in $G$ containing the $(u,v)$-link. By (\ref{spo}), the $(u,v)$-link is an edge, and the other two 
ends of $F$, say $a$ and $b$, belong to $X$. Moreover, since $X$ is maximal, $C-\{a,b\}$ 
has two components, each meeting $X$. Thus, $r\ge 4$, and assuming $a=x_{_0}$, we have 
$b=x_{_j}$, $1<j<r-1$. But then the union of $C$, the path $[x_{_1},v,x_{_{j+1}}]$, and 
$F$ is a non-\nsp P subdivision of $K_{_{3,3}}$, by \raf{k33}, contradiction.
\QED

\begin{lemma}\label{ear}
Let $G=(V,E)$ be virtually-$3$-connected and \nsp P. 
Let $J$ be an $S_{_r}$-subgraph of $G$, with the rim $C$, hub $v$, and $r=d^*(v)\ge 4$. 
Let $P$ be a $J$-ear with the ends $a,b\in C$. If each open $(a,b)$-segment of 
$C$ meets $N(v)$, then $r=4$, and either 

\mbox{(\sref{ear}.1)} $J\cup P\cong K_{_5}^-$, or 

\mbox{(\sref{ear}.2)} the sets $\{a,b\}$ and $N(v)$ are disjoint and nonadjacent in 
$C$, and $E(N(v))$ consists of \indent \indent \indent \ two disjoint edges, red with respect to $J\cup P$, lying
one in each component of $C-\{a,b\}$.
\end{lemma}

\noindent
{\sl Proof.}\\
(I) By Lemma~\sref{rd*}, the ends of $P$ should indeed belong to $C$. Put
$X\colon=N(v)=\{x_{_i}\colon\;i=0,\dots,r-1\;\mdl{r}\}$, in a circular order along $C$. 
Partition $C$ into semi-open $(a,b)$-segments $A$ (including $a$) and $B$ (including $b$), 
choosing the notation so as to maximize $|A\cap X|$. Clearly, $|A\cap X|\ge 2$. Let 
$x_{_0}$ belong to $B-b$, and $x_{_i}$ and $x_{_k}$ be the extreme members of $A\cap X$. 
Then, $(x_{_i}Ax_{_k})$ is an edge, by Lemma~\sref{xvy} (as applied to the path 
$[x_{_i},v,x_{_k}]$ and circuit $A\cup P$). Thus, $|A\cap X|=|B\cap X|=2$, whence 
$r=|X|=4$. Moreover, $|X\cap\{a,b\}|\ne 1$, by the maximality of $|A\cap X|$.\\

\noindent
(II) If $a,b\in X$, then the open segments $(x_{_{i-1}}Cx_{_i})$ are edges, for the same 
reason as in (I). Moreover, the length of $P$ is also $1$, for otherwise 
$J\cup P-\{av,bv\}$ is $\sbd{K_{_4}}{C_{_4}}$. Thus, $J\cup P\cong K_{_5}^-$.
and (\sref{ear}.1) is satisfied.\\

\noindent
(III) If $a,b\not\in X$, assume $a\in(x_{_3}Cx_{_0})$ and $b\in(x_{_1}Cx_{_2})$; again, 
Lemma~\sref{xvy} implies that $(x_{_0}Cx_{_1})$ and $(x_{_2}Cx_{_3})$ are edges. To show 
that no other edge is spanned by $N(v)$, return to (I) and (II) assuming that $P$ is an edge spanned by $N(v)$.\\

\noindent
(IV) Suppose finally that $\{a,b\}$ and $N(v)$ are adjacent in $C$, say $(aCx_{_0})$ is
an edge, by symmetry. This contradicts $J\cup P-\{vx_{_1},x_{_2}x_{_3}\}$ being 
\nsp P. Thus, (III) and (IV) imply (\sref{ear}.2). \QED

\begin{lemma}\label{tri}
Let $G=(V,E)$ be virtually $3$-connected and \nsp P.
A vertex $v\in V^*$ with $d^*(v)>1$ belongs to a triangle, unless $G$ is $S_{_r}$,
$r\ge 3$, with skeleton $K_{_{1,r}}$.
\end{lemma}

\noindent
{\sl Proof.}
For $d^*(v)=2$, the assertion follows from (\sref{rd*}.2); so suppose, to the contrary, 
that $d^*(v)\ge 3$ and $v$ belongs to no triangle.

Let $J$ be an $S_{_r}$-subgraph of $G$ with the hub $v$ and $r=d^*(v)$ with $N(v)$ on its rim as guaranteed by Lemma~\ref{rd*}. We show that 
the only skeleton vertices in $C\colon=J-v$ are the members of $N(v)$. If not, choose 
$a,b\in N(v)$ such that the open segment $(aCb)$ contains no member of $N(v)$ and 
meets $V^*$. Since $G$ is virtually $3$-connected, there is a $J$-ear $P$ linking 
the components of $C-\{a,b\}$, with the ends $s\in (aCb)$ and $t$. Since 
$|N(v)|=d^*(v)\ge 3$, at least one closed $(s,t)$-segment of $C$ contains more than one 
member of $N(v)$. If $Q$ is such a segment, and $x,y\in N(v)\cap Q$, then the path 
$[x,v,y]$ and the circuit $Q\cup P$ satisfy \raf{DN}  (because $x$ and $y$ are 
nonadjacent), so that $G$ is non-\nsp P, contradiction.
  
It remains to show that $J$ admits no ear. Indeed, if $P$ is an $(s,t)$-ear of $J$, then 
$s,t\in N(v)$. Since $G$ is simple, both open $(s,t)$-segments of $C$ meet $N(v)$, and 
we have the same contradiction as above. Thus, $G$ is $S_{_r}$ and $G^*\cong K_{_{1,r}}$.
\QED 
\\

Except for quite special cases, no edge of $G$ belongs to two skeleton triangles to
form $K_{_4}^-$:

\begin{lemma}\label{k4-}
Let $G=(V,E)$ be virtually-$3$-connectes and \nsp P. 
If $G^*$ contains $K_{_4}^-$, then either $G$ contains $K_{_5}^-$ or 
$S_{_r}\subseteq G\subseteq\modi{W_{_r}}$ with $r\ge 4$.
\end{lemma}

\noindent
{\sl Proof.}
Let $J$ be a $K_{_4}^-$-subgraph of $G^*$ with $2$-valent vertices $x,y$ and
$3$-valent vertices $u,v$, and suppose that $d^*(u)\le d^*(v)$. 

We show that since $|V^*|>4$, we have, by symmetry, that $d^*(v)\ge 4$. To see this, assume, to the contrary, that $d^*(u) = d^*(v) = 3$, and let $H$ be an $S_{_r}$-subgraph, $r=d^*(v)$, of $G$ with the hub $v$ and rim $C$ containing $N^2(v) \subseteq N^1(v)$. We may choose $H$ such that $C$ contains the edges $\{ux,uy\}$. If the open $(x,y)$-segment of $C$ not containing $u$, namely $L$, contains a skeleton vertex, say, $t$, then $d^*(v) \geq 4$ or $d^*(u) \geq 4$. Indeed, a $(t,J)$-fan $F$ with ends in $x,y$ also ends in one of $\{u,v\}$. Clearly, the $tFz$, $z \in \{u,v\}$, member of $F$ is an edge or a forbidden topological $K_{_4}$ appears in $J \cup F - uv$; implying that $d^*\geq 4$ for at least one of $\{u,v\}$. 
 
Thus, we may assume that $L$ is a link. Since $|V^*| > 4$, there is a skeleton vertex $t$ external to $J$. A $(t,J)$-fan $F'$ in this case must end at $\{u,v,w\}$, where $w \in \{x,y\}$. Otherwise, we reach the same contradiction as before. 
The assumption that $d^*(v) = d^*(u) = 3$ implies that the $uF't$ and $vF't$ members of $F'$ are not edges. In addition, $wF't$ is an edge or a forbidden topological $K_{_4}$ appears. Thus, $J \cup F' \cup L - uw' - vw$ is a forbidden topological $K_{_4}$, where $\{w'\} = \{x,y\} \sm \{w\}$.

We have seen that at least one of $u,v$ has $d^* \geq 4$, we assume it is $v$. 

Let $H$ be as above and let $H'$ be $\modi{W_{_r}}$ containing $H$. Then, either $G\subseteq H'$, or $H$ has an ear as assumed in Lemma~\sref{ear}. In the latter case, the lemma implies that $G$ contains $K_{_5}^-$. 
\QED

\subsection{Structure surrounding a skeleton triangle}\label{triangle} Throughout this section 
\begin{equation}\label{T}
\mbox{$T\subseteq E^*$ denotes a skeleton triangle, with $V(T)=\{x,y,z\}$.}
\end{equation}

\begin{lemma}\label{thre}
Let $G=(V,E)$ be virtually-$3$-connected and \nsp P and let $T$ be as in \raf{T}.
At least two edges of $T$ are red. Also, if one edge of $T$ has a parallel thread attached to its ends, then $S_{_r}\subseteq G\subseteq\modi{W_{_r}}$, with $r\ge 3$.
\end{lemma}

\noindent
{\sl Proof.}
Suppose $P$ is an $(x,y)$-thread. Since $|V^*|>4$, 
the subgraph $G-z-xy$ is $2$-connected, so that there exists a circuit $C$ in $G-z$ 
containing $P$. By Lemma~\sref{xvy}, $N^2(z)\sm\{x,y\}=\empt$ and $N^1(z)\subseteq V(C)$ so that $d^*(z)\geq 3$ (or $\{x,y\}$ disconnect $z$ contradicting virtual $3$-cnnectivity). 
Denote by $J$ the subgraph composed of $C\cup T$ and the edges incident with $z$, so 
that $r=d^*(z)$; we show that $G$ is isomorphic to a subgraph of $\modi{W_{_r}}$. If not, then $G$ contains a $J$-ear $Q$, with some ends $a,b$; clearly $a,b\in V(C)$. Since 
$G$ is virtually $3$-connected, $Q$ may be chosen so that each open $(a,b)$-segment of 
$C$ meets $N(z)$. Then, $C\cup Q$ contains a circuit $C'$ which together with the path 
$[x,z,y]$ satisfies \raf{DN}; so that $G$ is non-\nsp P, by Lemma~\sref{xvy}, a contradiction.
Thus, $S_{_r}\subseteq G\subseteq\modi{W_{_r}}$, with the hub $z$. Since $|V^*|>4$, we 
have $d^*(z)\ge 4$. Then, the edges $xz$ and $yz$ are red, by \raf{spo}.

It remains two show that at least two edges of $T$ are red. We may assume that at least one edge of $T$ is not red, otherwise the assertion follows. Above we have seen that attaching a thread in parallel to one of these edges implies that $z$ is the hub of an $S_{_r}$-subgraph where $r = d^*(z) \geq 4$, and that the added thread lies on the rim of the $S_{_r}$ subgraph. We may replace this thread by the edge of $T$ assumed non-red and obtain that, in $G$, $z$ is a hub of an $S_{_r}$-subgraph where $r = d^*(z) \geq 4$. The claim then follows by \raf{spo}.\QED
 
\begin{lemma}\label{tri1}
Let $G=(V,E)$ be virtually-$3$-connected and \nsp P and let $T$ be as in \raf{T}.
Let $C$ be a circuit in $G-xy$ traversing $x$ and $y$. If $z$ is not in $C$, then 
$G$ contains $K_{_5}^-$.
\end{lemma}

\noindent
{\sl Proof.}
Suppose that $z\not\in V(C)$; then $N^2(z)\sm\{x,y\}=\empt$ and $N^1(z)\subseteq V(C)$, 
by Lemmas~\sref{xvy} and~\sref{thre}. Then, $d^*(z)\ge 3$, so that there exists an  
$s\in N^1(z)\cap C$ distinct from $x$ and $y$. Let $P$ and $Q$ denote the open 
$(x,y)$-segments of $C$, and suppose that $s$ lies in $P$. Since $G$ is virtually
$3$-connected and $Q$ meets $V^*$ (by Lemma~\sref{thre}), $C\cup T$ has an ear
linking some vertices $t$ of $Q$ and $u$ of $P\cup\{v\}$. Let $J$ denote the subgraph
composed of $C\cup T$, this ear, and the edge $zs$. We show that $J\cong K_{_5}^-$.

First, $u=z$. Indeed, otherwise we may assume $z\in(xPs]$, by symmetry. Then, $(sPy)$ is 
an edge, by Lemma~\sref{xvy}. The same lemma, applied to the path $[s,y,x]$ and circuit 
$[x,z,s]\cup P$, implies that $u=s$ and $P$ is an edge; the latter contradicts the choice 
of $C$. 

Second, Lemma~\sref{xvy} may now be applied to $J-(tQx)-(y,z)$ and $J-(tQy)-(x,z)$, and 
we conclude that $(sPx)$ and $(sPy)$ are edges. Finally, the paths of $J$ linking $t$
to $T$ are edges: by symmetry, it suffices to observe that if $(tQx)$ is not an edge,
then $J-\{sx,yz\}$ is non-\nsp P. Thus $J\cong K_{_5}^-$, as required.
\QED

\begin{lemma}\label{tri2}
Let $G=(V,E)$ be virtually-$3$-connected and \nsp P and let $T$ be as in \raf{T}.
Suppose that $|N(x)|\le|N(y)|\le|N(z)|$. Then, either

\mbox{(\sref{tri2}.1)} $G$ has a $K_{_5}^-$-subgraph, or

\mbox{(\sref{tri2}.2)} $S_{_r}\subseteq G\subseteq\modi{W_{_r}}$, $r\ge 3$, or

\mbox{(\sref{tri2}.3)}
$|N(x)|=|N(y)|=3$, and 

\mbox{(\sref{tri2}.4)}
if $|N(z)|>3$, then $d^*(x)=d^*(y)=2$ and $d^*(z)=2$ or $4$.
\end{lemma}

\noindent
{\sl Proof.}
Suppose that $G$ satisfies neither of (\sref{tri2}.1-2). To show (\sref{tri2}.3), note
first that 

\begin{equation}\label{gtc}
\mbox{\rm if $|N(y)|>3$, then $G-T$ has a circuit containing $V(T)$.} 
\end{equation}

\noindent
{\sl Proof.}
Suppose \raf{gtc} is false, then the subgraph $G-T$ is not $2$-connected. Indeed, if $G-T$ has no $V(T)$-circuit
and is $2$-connected, then it contains a subdivision $J$ of $K_{_{3,2}}$ with $V(T)$ as 
the larger part, by Observation~\sref{whi}. Lemma~\sref{xvy} then implies that $J\cup T\cong K^-_{_5}$, so that $G$ satisfies (\sref{tri2}.1). Hence, $G-T=H'\cup H''$, with $H'\cap H''$ consisting of some 
$a\in V^*\sm V(T)$, and we assume $|V(T)\cap H''|=1$. 

Since $G$ is virtually $3$-connected, $(H''-a)\cap V^*$ consists of just one vertex which
belongs to $T$ and is $3$-valent. If $|N(y)|>3$, this may only be $x$, 
so that $|N(x)|=3$. Then, $G-x-yz$ is $2$-connected, for otherwise 
$G-x-yz=L'\cup L''$ with $L'\cap L''$ consisting of one vertex, say $b$, and we may 
assume that $z\in L'-b$, $y\in L''-b$, and $a\in L''$. Since $\{b,y\}$ does not 
disconnect $G$, we have $(L''-b)\cap V^*=\{y\}$ whence $|N(y)|=3$, contradiction. 

Thus, $G-x-yz$ has a circuit $C$ traversing $y$ and $z$, so that $G$ contains 
$K_{_5}^-$, by Lemma~\sref{tri1}; this contradiction proves \raf{gtc}.\inQED \\

Let now $C$ be a circuit of $G-T$ containing $V(T)$, and put $J\colon=C\cup T$; suppose 
$(yCz)\cap V^*\ne\empt$, by Lemma~\sref{thre}. Since $G$ is virtually $3$-connected, 
it contains an ear $P$ of $J$ linking the components of $J-\{y,z\}$; let the ends of $P$ 
be $s\in(yCz)$ and some $t$. There are two possibilities.

First, $t\ne x$, so that, say, $t\in (xCz)$. Then, $(sCz)$ and $(tCz)$ are edges, by 
Lemma~\sref{xvy}, so that $d^*(z)\ge 4$. By Lemma~\sref{ear}, there is no ear of 
$J\cup P$ with the ends in $C'\colon=C\cup P-\{sz,tz\}$, so that 
$S_{_r}\subseteq G\subseteq\modi{W_{_r}}$ with $r=d^*(z)$, contradiction.

Second, $t=x$. Then, the segments $(sCy)$ and $(tCy)$ are edges, the latter by  
Lemma~\sref{xvy} as applied to the path $[x,z,y]$ and the circuit of $C\cup P-z$, and 
the former by symmetry. Thus, $G$ contains $K_{_4}^-$, and therefore satisfies one of 
(\sref{tri2}.1-2), by Lemma~\sref{k4-}, contradiction. Thus, $|N(x)|=|N(y)|=3$. 

To show (\sref{tri2}.4), let $|N(z)|\ge 4$; it suffices to show that $d^*(x)=d^*(y)=2$. 
Suppose, to the contrary, that $d^*(x)=3$. Let $s\in N^1(x)\sm V(T)$, and let 
$t\in N(y)\sm V(T)$, with the $(t,y)$-link $P$. Since $G$ is virtually $3$-connected, 
$s\ne t$. We show that $G-\{x,y\}$ has a circuit traversing $s$ and $z$, so that $G$ 
is non-\nsp P, by Lemma~\sref{xvy}. Indeed, if no such circuit exists, we have 
$G-\{x,y\}=H'\cup H''$ where $H'\cap H''$ consists of one vertex, say $u$, while $s$ 
and $z$ belong to $H'-u$ and $H''-u$ respectively. Since $|N(z)\sm\{x,y\}|\ge 2$, 
the pair $u$, $z$ disconnects $G$, contradiction. 
\QED

\section{Proof of Theorem~\sref{B}}\label{proof-B}

\noindent
{\sc Only if.}
Let $G$ be an ASP graph with $|V^*|>6$. If $G$ has a $\sbd{K_{_4}}{C_{_4}}$-subgraph,
then it satisfies (\sref{B}.1), by Proposition~\sref{kc4}. Suppose now that $G$ is \nsp P.\\

\noindent
(I) Since $|V^*|>5$, $G$ contains no $K_{_5}^-$-subgraph, because a fan with its ends in such a subgraph yields a non-\nsp P graph.\\

\noindent
(II) If $\gD(G^*)\ge 5$, then $G$ satisfies (\sref{B}.2), by Lemmas~\sref{rd*} and 
\sref{ear}.\\

\noindent
(III) Suppose then that $\gD(G^*)\le 4$, $G$ contains no $K_{_5}^-$-subgraph, and satisfies neither of (\sref{B}.1-2). We show that $G$ is a fishpond. This is clearly true if $d^*(v) \leq 1$ for each $v \in V$. So, let $v\in V^*$ have $d^*(v)\ge 2$. Then, (A.1) follows from Lemma~\sref{tri}, and (A.2) from Lemmas~\sref{thre} and  
\sref{tri2}. To show (A.3), let now $J$ be a component of $G^*$ with $|V(J)|>3$. By (II), $\gD(J)\leq 4$, by (A.1), which we know to be true now, $J$ contains a triangle, and it being connected and $|V(J)| > 3$ implies that $\Delta(J) \geq 3$. 
Thus, it remains to consider to case $\Delta(J) = 3$ or $4$. 

Suppose that $\gD(J)=3$, and let $x\in V(J)$ have $d^*(x)\ge 2$. By Lemma~\sref{tri}, $x$ belongs to a triangle of $J$, whose vertices are, say, $x,y,z$. Since
$\gD(J)=3$, one of these vertices has $d^*=3$. Then, $|N(x)|=|N(y)|=|N(z)|=3$, by
Lemma~\sref{tri2}, so that $G$ satisfies (A.3.2). 

Suppose finally that $\gD(J)=4$, and choose $v\in V(J)$ with $d^*(v)=4$. By 
(\sref{rd*}.1), $J$ contains the skeleton of an $S_{_4}$-subgraph $H$ of $G$ with hub 
$v$. Since $G\not\subseteq\modi{W_{_4}}$ and $K_{_5}^-\not\subseteq G$, there exists 
an ear $P$ of $H$ with the ends in the rim. By Lemma~\sref{ear}, the skeleton of $H$ is a union of two edge-disjoint triangles containing $v$, as in (A.3.3), and coincides with $J$. 

Thus $J$ is a fishpond, as required.\\

\noindent
{\sc If.} Cases (\sref{B}.1-2) are trivial; it remains to show that a fishpond is \nsp P. This is clearly so if $d^*(v) \leq 1$ for each vertex $v$. So, let $G$ be a fishpond, and $[x,v,y]$ be a path of length $2$ in $G^*$. By Lemma~\sref{xvy}, it suffices to show that when $G-v$ 
has a pair of internally disjoint $xy$-paths of length $>1$, their union contains 
$N^1(v)$, and $N^2(v)\subseteq\{x,y\}$.  

By (A.1), $v$ belongs to a triangle, say $T$. If $x$ and $y$ are vertices of $T$, then
we may assume that $|N(x)|=3$, by (A.2), so that $G-v-xy$ has no circuit traversing 
$x$, and the requirement is trivially satisfied. Suppose now that $x$ and $y$ are 
nonadjacent, so that the vertices of $T$ are $v$, $x$ and some $z\ne y$. Then, 
$d^*(v)\ge 3$. By (A.2), we have $N^2(v)\subseteq\{y\}$. By (A.3.3), if $d^*(v)=4$, 
then $v$ is a common vertex of two edge-disjoint triangles. Then, (A.2) implies that 
$v$ is regular, and any circuit in $G-v$ traversing $x$ and $y$ contains $N^1(v)$, as 
required. \QED

\section{Colorability of ASP and \nsp P graphs}\label{color-ASP}
A receptacle is called {\sl extreme} if it has exactly one window.
Every graph $G$ with $\kappa(G)=2$ and $\delta(G)\geq 3$ has at least two extreme receptacles.
If $H$ is an extreme receptacle of $G$ with $\{u,v\}$ as its sole window, then $\kappa(H+uv) \geq 3$ (if $uv\in E(H)$, then $\kappa(H) \geq 3$). The union of $H$ with a $(u,v)$-thread internally disjoint of $H$ is a virtually $3$-connected graph with precisely one thread.

\subsection{Colorability of ASP graphs}\\

\noindent
\bfm{Proof of Theorem~\sref{C}.} Any graph of order $\leq 6$ is contained in $K^-_{_6}$ and is $5$-colorable. It remains to consider ASP graphs of order $> 6$. Let, to the contrary, $G$ be a minimal non-$5$-colorable (i.e., $6$-critical) ASP graph of such an order. Then, $\delta(G)\ge 5$, $\kappa(G) \geq 2$, and any disconnector of $G$ does not induce a complete graph.

By Theorem~\sref{A}, $\kappa(G)=2$. Let $H$ be an extreme receptacle of $G$ with the sole window $\{x,y\}$. Let $H'$ be the union of $H$ and an $xy$-path in $G$ internally disjoint of $H$. Such a path is not an edge, implying that $H'$ is a virtually-$3$-connected ASP graph with precisely one thread. 

Since $\delta(G) \geq 5$, each vertex in $V(H) \sm \{x,y\}$ is at least $5$-valent. By Theorems~\ref{A} and~\ref{B}, the skeleton of $H'$  has at most $6$ vertices, for all other options have vertices at most $4$-valent. On the other hand, $\delta(G) \geq 5$ implies that $|V(H)| \geq 6$. Hence, $H$ is isomorphic to $K^-_{_6}$ so that $H'$ is a $K^-_{_6}$ with a thread attached to the nonadjacent vertices. Such a graph is non-ASP.    

It is quite standard now to see that the above proof yields a polynomial time algorithm for $5$-vertex-coloring ASP graphs using techniques presented in~\cite{bm, diestel}.\QED \\

\noindent
\scaps{$5$-chromatic ASP graphs.} We present now a construction of $5$-chromatic ASP graphs containing no $K_{_5}$ as a subgraph: Let $J$ be an arbitrary $2$-connected graph. Replace each edge $e$ of $J$, except for one edge $e'$, with a copy of $K_{_5}^-$ by identifying the ends of $e$ with the pair of nonadjacent vertices in $K_{_5}^-$. Every receptacle of the resulting graph $J'$ is extreme and isomorphic to $K^-_{_5}$. As $K^-_{_5}$ with a thread linking its nonadjacent vertices appended is ASP, $J'$ is ASP by Theorems~\sref{A}, \ref{B} and \raf{rcpt}. If $e'$ belongs to a circuit of $J$, then the resulting graph is not $4$-colorable. This follows from the fact that every $4$-coloring of $K_{_5}^-$ has the nonadjacent vertices colored the same.

The above construction leads to a proof of Theorem~\ref{D}.\\

\noindent
\bfm{Proof of Theorem~\sref{D}.}
Let $X_1$ and $X_2$ be two vertex-disjoint copies of $K^-_{_5}$, and let $x_i$, $y_i$ denote the ends of the removed edge in $X_i$, $i=1,2$. Define $Y$ to be the graph obtained from $X_1$ and $X_2$ by adding the edge $x_1x_2$. In every $4$-coloring of $Y$, $y_1$ and $y_2$ are colored distinctively. Also, $Y$ and a thread linking $y_{_1}$ and $y_{_2}$ appended is ASP, by Theorems~\ref{A}, \ref{B}, and \raf{rcpt}. 

Let $J$ be an arbitrary $2$-connected graph. Replace each edge $e$ of $J$ with a copy 
of $Y$, by identifying its $\{y_1,y_2\}$-pair with the ends of $e$. Since any receptacle of the resulting graph $J'$ is either isomorphic to $Y$ or to $K^-_{_5}$ (the latter if the receptacle is extreme), then $J'$ is ASP, by Theorems~\ref{A}, \ref{B}, and \raf{rcpt}. It follows that $\chi(J) \leq 4$ if and only if $\chi(J') \leq 4$. Checking this inequality for an arbitrary $2$-connected graph is known to be NP-complete~\cite{gj}. \QED
 
\subsection{Colorability of \nsp P graphs} In a similar manner to \raf{rcpt}, we have that
\begin{equation}\label{rcpt2}
\mbox{\rm a graph is \nsp P if and only if each of its receptacles is \nsp P.}
\end{equation}

\noindent
\bfm{Proof of Proposition~\sref{E}.} Any graph of order $\leq 5$ contained in $K^-_{_5}$ is $4$-colorable. It remains to consider ASP graphs of order $> 5$. Let, to the contrary, $G$ be a minimal non-$4$-colorable (i.e., $5$-critical) \nsp P graph of such an order. Then, $\delta(G)\ge 4$, $\kappa(G) \geq 2$, and any disconnector of $G$ does not induce a complete graph. 

By Theorem~\ref{A}, $\kappa(G)=2$. Let $H,H'$, and $\{x,y\}$ be as in the proof of Theorem~\ref{C}.
Since $\delta(G) \geq 4$, each vertex in $V(H) \sm \{x,y\}$ is at least $4$-valent. By Theorems~\ref{A} and~\ref{B}, the skeleton of $H'$ has at most $6$ vertices, for all other options have vertices at most $3$-valent. On the other hand $|V(H)| \geq 5$, as $\delta(G) \geq 4$. So that the skeleton of $H'$ contains $K^-_{_4}$. Lemma~\ref{k4-} and \raf{rcpt2} then imply that $H'$ is $W_{_r}$, $r \geq 4$, with precisely one of the rim edges replaced by a thread; so that $V(H) \sm \{x,y\}$ contains $3$-valent vertices; a contradiction. 

As in Theorem~\ref{D}, this proof outlines a polynomial time algorithm for $4$-coloring \nsp P graphs. \QED \ \\

\noindent
\scaps{$4$-chromatic \nsp P graphs.} We present now a construction of planar $4$-chromatic \nsp P graphs containing no $K_{_4}$ as a subgraph. 
For $r \geq 3$, let $W^-_{_r}$ denote $W_{_r}$ with a single rim edge removed.
Let $J$ be an arbitrary planar $2$-connected graph. Replace each edge $e$ of $J$, except for one edge $e'$, with a copy of $W^-_{_r}$, $r \geq 3$ and odd, by identifying the ends of $e$ with the ends of the removed edge in $W^-_{_r}$. The resulting graph $J'$ is clearly planar and every receptacle of $J'$ is extreme and isomorphic to $W^-_{_r}$. As $W^-_{_r}$ with a thread linking the ends of its removed rim edge is \nsp P, $J'$ is \nsp P by Theorems~\sref{A},\ref{B} and \raf{rcpt2}. If $e'$ belongs to a circuit of $J$, then the resulting graph is not $4$-colorable. This follows from the fact that every $3$-coloring of $W^-_{_r}$ has the ends of its removed edge colored the same.

The above construction leads to a proof of the following. \\

\noindent
\bfm{Proof of Proposition~\sref{F}.} For $r \geq 4$ and even, any $3$-coloring of $W_{_r}^-$ has the ends of the removed edge colored distinctively. Let $J$ be an arbitrary $2$-connected planar graph. Replace each edge $e$ of $J$ with a copy of $W_{_r}^-$, by identifying the ends of its removed edge with the ends of $e$. The resulting 
graph is planar and \nsp P, by Theorems~\ref{A}, \ref{B}, and \raf{rcpt2}. The resulting graph is $3$-colorable if and only if $\chi(J)\le 3$. Checking this inequality for an arbitrary $2$-connected planar graph is known to be NP-complete \cite{gj}.\QED \\

\end{document}